\documentclass [11pt,oneside]{article}

\reversemarginpar \textwidth=14cm \textheight=19cm

\usepackage{amssymb} \usepackage{amsfonts} \usepackage{amsmath}
\usepackage{amsthm} \usepackage{epsfig}

\newcommand{\mettifig}[1]{\epsfig{file=#1}}

\newtheorem{lemma}{Lemma}[section]
\newtheorem{teo}[lemma]{Theorem}
\newtheorem{rem}[lemma]{Remark} 
\newtheorem{prop}[lemma]{Proposition}

\newcommand{\matR} {\ensuremath {\mathbb{R}}}
\newcommand{\matQ} {\ensuremath {\mathbb{Q}}}
\newcommand{\matZ} {\ensuremath {\mathbb{Z}}}
\newcommand{\matC} {\ensuremath {\mathbb{C}}}
\newcommand{\matP} {\ensuremath {\mathbb{P}}}
\newcommand{\matH} {\ensuremath {\mathbb{H}}}
\newcommand{\matRP} {\ensuremath {\mathbb{RP}}}

\newcommand{\calR} {\ensuremath {\mathcal{R}}}

\def\closure#1{{\rm Cl}(#1)}

\newcommand{\ptwoirred}{$\matP^2$-irreducible}

\newcommand{\timtil}{\begin{picture}(12,12)
\put(2,0){$\times$}\put(2,4.5){$\sim$}\end{picture}}

\newcommand{\matr} [4] {\tiny{\left(\begin{array}{@{}c@{\ }c@{}} #1 & #2 \\ #3 & #4 \\ \end{array} \right)}}

\usepackage{caption}

\author{Gennaro \textsc{Amendola} \and Bruno \textsc{Martelli}}

\title{Non-orientable $3$-manifolds of small complexity}

\begin{document}

\maketitle

\begin{abstract}
        \noindent
  We classify all closed non-orientable \ptwoirred\ $3$-manifolds having
  complexity up to $6$ and we describe some having complexity $7$. 
  We show in particular that there is no such manifold
  with complexity less than $6$, and that those having complexity $6$
  are precisely the $4$ flat non-orientable ones and the filling of
  the Gieseking manifold, which is of type Sol.
  The manifolds having complexity $7$
  we describe are Seifert manifolds of type $\matH^2\times S^1$ and a manifold of type
  Sol.\\[4pt]
\end{abstract}

\section{Introduction} \label{complexity:section}
In~\cite{Mat} Matveev defined
for any compact 3-manifold $M$ a non-negative integer $c(M)$, 
which he called the \emph{complexity} of $M$. The complexity function $c$
has the following remarkable properties: it is additive on connected sums,
it does not increase when cutting along incompressible surfaces, and it 
is finite-to-one on the most interesting sets of 3-manifolds. Namely,
among the compact 3-manifolds having complexity $c$ 
there is only a finite number of 
closed \ptwoirred\ ones,
and a finite number of hyperbolic ones 
(with cusps and/or with geodesic boundary). 
At present, hyperbolic
manifolds with cusps are classified in~\cite{CaHiWe} for $c\leqslant
7$, and orientable hyperbolic manifolds with
geodesic boundary are classified in~\cite{FriMaPe2} for
$c\leqslant 4$.
In this paper we concentrate on the closed \ptwoirred\ case: 
the complexity of such an $M$ 
is then precisely the minimal number of tetrahedra needed to triangulate
it, except when $c(M)=0$, \emph{i.e.}~when $M$ is $S^3,\matRP^3$ or $L_{3,1}$.

\paragraph{Known results on closed manifolds}
We recall that there are 8 important 3-dimensional geometries, six of them 
concerning Seifert manifolds. 
The geometry of a Seifert manifold is determined by two invariants of any of its fibrations,
namely the Euler characteristic $\chi^{\rm orb}$ of the base orbifold and the Euler number
$e$ of the fibration, according to Table~\ref{tabellina}. The two non-Seifert
geometries are the hyperbolic and the Sol ones.

\begin{table}
\begin{center}
\begin{tabular}{c|ccc} 
\phantom{\Big|} & $\chi^{\rm orb}>0$ & $\chi^{\rm orb}=0$ & $\chi^{\rm orb}<0$ \\ \hline
\phantom{\Big|} $e=0$ & $S^2\times\matR$ & $E^3$ & $H^2\times\matR$ \\
\phantom{\Big|} $e\neq 0$ & $S^3$ & Nil & $\widetilde{{\rm SL}_2\matR}$ \\ 
\end{tabular}
\end{center}
\caption{The six Seifert geometries.}
\label{tabellina}
\end{table} 

Using computers, closed \emph{orientable} irreducible 3-manifolds having 
complexity up to $6$~\cite{Mat} and then up to $9$~\cite{MaPe} have been classified.
The complete list is available from~\cite{weblist}, and we summarize it in 
the first half of Table~\ref{tabella}.
In particular, the orientable manifolds with
$c\leqslant 5$ are Seifert with $\chi^{\rm orb}>0$, and 
those with $c\leqslant 6$ are Seifert with $\chi^{\rm orb}\geqslant 0$,
including all $6$ flat ones.
Seifert manifolds with $\chi^{\rm orb}<0$ or Sol geometry appear with $c=7$, and 
the first hyperbolic ones have $c=9$ (this was first proved
in~\cite{MaFo}). Manifolds with non-trivial JSJ decomposition
appear with $c=7$: 
each such manifold with $c\leqslant 9$ actually decomposes into Seifert pieces.
We show in Section~\ref{teoria:section} that the first manifold
whose JSJ decomposition is non-trivial and contains a hyperbolic piece
has $c\leqslant 11$, and we explain why we think it should have $c=11$. 

\begin{rem}
{\em The number of manifolds having complexity $9$ is 1155. The wrong number 1156
found in~\cite{MaPe} was the result of a list containing the same graph manifold twice.
}
\end{rem}

\paragraph{Main statement}
We prove in this paper that non-orientable \ptwoirred\ manifolds, more
or less, follow the same scheme.
Taking into account that a non-orientable Seifert manifold has
Euler number zero~\cite{Sco}, we mean the following.
\begin{teo} \label{main:teo}
\begin{itemize}
\item
  There are no closed non-orientable \ptwoirred\
  manifolds with $c\leqslant 5$,
\item
  the only ones with $c=6$ are the $4$ flat ones and the torus bundle
  over $S^1$ of type {\rm Sol} with monodromy trace $\matr 1110$,
\item
  there are some of type $\matH^2\times\matR$ and of type {\rm Sol} with $c=7$.
\end{itemize}
\end{teo}

These results are summarized in the second half of Table~\ref{tabella}.
We emphasize that the proof of Theorem~\ref{main:teo} 
is theoretical (\emph{i.e.}~it makes no use of any computer result).
We end this
section by defining Matveev's complexity and by describing the main line of
the proof. Some techniques taken 
from~\cite{MaPe} will be briefly summarized in Section~\ref{teoria:section},
and these techniques will be used in Section~\ref{fine_dimo:sec} to
conclude the proof. The proofs of some technical lemmas are postponed to 
Section~\ref{proof_lemmas:sec}.

\begin{table}[t]
  \begin{center}
    \begin{tabular}{r|c|c|c|c|c|c|c|c|c|c}
      & $0$ & $1$ & $2$ & $3$ & $4$ & $5$ & $6$ & $7$ & $8$ & $9$ \\
      
      \hline\hline

      \multicolumn{11}{c}{orientable\phantom{\Big|}} \\

      \hline
      \hline
      
      \phantom{\Big|} 
      lens &
      $3$ & $2$ & $3$ & $6$ & $10$ & $20$ & $36$ & $72$ & $136$ & $272$ \\

      \hline

      other elliptic &
      \multicolumn{1}{c}{\phantom{\Big|}} & & $1$ & $1$ & $4$ & $11$ & $25$ & $45$ & $78$ & $142$ \\

      \hline
      
      flat &
      \multicolumn{5}{c}{\phantom{\Big|}} & & $6$ & \multicolumn{3}{c}{} \\
      
      \hline
      
      Nil &
      \multicolumn{5}{c}{\phantom{\Big|}} & & $7$ & $10$ & $14$ & $15$ \\

      \hline

      $\widetilde{{\rm SL}_2(\matR)}$ &
      \multicolumn{6}{c}{\phantom{$\Big|^|$}} & & $39$ & $162$ & $514$ \\
      
      \hline

      Sol &
      \multicolumn{6}{c}{\phantom{\Big|}} & & $5$ & $9$ & $23$ \\

      \hline

      $\matH^2\times\matR$ &
      \multicolumn{7}{c}{\phantom{\Big|}} & & $2$ & \\
      
      \hline 

      hyperbolic &
      \multicolumn{8}{c}{\phantom{\Big|}} & & $4$ \\
      
      \hline

      non-trivial JSJ &
      \multicolumn{6}{c}{\phantom{\Big|}} & & $4$ & $35$ & $185$ \\
      
      \hline
      \hline 

      total orientable \phantom{\Big|} &
      $3$ & $2$ & $4$ & $7$ & $14$ & $31$ & $74$ & $175$ & $436$ & $1155$ \\

      \hline
      \hline



      \multicolumn{11}{c}{{\bf non-orientable}\phantom{\Big|}} \\

      \hline
      \hline

      {\bf flat} &
      \multicolumn{5}{c}{\phantom{\Big|}} & & $\bf{4}$ & \multicolumn{3}{c}{} \\

      \hline

      {\bf Sol} &
      \multicolumn{5}{c}{\phantom{\Big|}} & & $\bf{1}$ & $\bf{>0}$ & \multicolumn{2}{c}{?} \\

      \hline

      {\bf $\matH^2\times\matR$} &
      \multicolumn{6}{c}{\phantom{\Big|}} & & $\bf{>0}$ & \multicolumn{2}{c}{?} \\ 

      \hline
      \hline 

      {\bf total non-orientable} &
      \multicolumn{5}{c}{\phantom{\Big|}} & & $\bf{5}$ & $\bf{>0}$ & \multicolumn{2}{c}{?} \\

      \hline

    \end{tabular}
  \end{center} 
  \caption{The number of \ptwoirred\ manifolds of given complexity (up to $9$) and geometry
    (empty boxes contain $0$).}
  \label{tabella}
\end{table}

\begin{figure}
\begin{center}
\mettifig{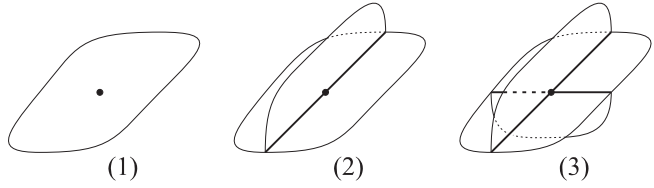}
\caption{Neighborhoods of points in a standard polyhedron.} 
\label{standard_nhbds:fig}
\end{center}
\end{figure}

\paragraph{Definition of complexity}
A compact 2-dimensional
polyhedron $P$ is said to be \emph{simple} if the link of every point in $P$ is contained in
the 1-skeleton $K$ of the tetrahedron. A point, a compact graph, a
compact surface are thus simple. Three important possible kinds of
neighborhoods of points are shown in Fig.~\ref{standard_nhbds:fig}. 
A point having the whole of $K$ as a link is 
called a \emph{vertex}, and its regular neighborhood
is shown in Fig.~\ref{standard_nhbds:fig}-(3). 
The set $V(P)$ of the vertices of $P$ consists of isolated points, so
it is finite. Points, graphs and surfaces of course do not contain vertices.
A compact polyhedron $P\subset M$ is a \emph{spine} of the closed manifold $M$
if $M\setminus P$ is an open ball. The \emph{complexity} $c(M)$ of
a closed 3-manifold $M$ is then defined as the minimal number of vertices of
a simple spine of $M$.

Now a point is a spine of $S^3$, the projective plane $\matRP^2$ is a
spine of $\matRP^3$ and the ``triple hat'' -- a triangle with all edges
identified in the same direction -- is a simple spine of
$L_{3,1}$. Since these spines do not contain vertices, we have 
$c(S^3)=c(\matRP^3)=c(L_{3,1})=0$. In general, to calculate the
complexity of a manifold we must look for its \emph{minimal} spines,
\emph{i.e.}~the simple spines with the lowest number of vertices.
It turns out~\cite{Mat, MaPe:nonori} that if $M$ is \ptwoirred\ and
distinct from $S^3, \matRP^3, L_{3,1}$ then it has a minimal spine
which is \emph{standard}. A polyhedron is standard when
every point has a neighborhood of one of the types (1)-(3) shown in
Fig.~\ref{standard_nhbds:fig}, and the sets of such points induce a
cellularization of $P$. That is, defining $S(P)$ as the set of points of type (2) or (3),
the components of $P \setminus S(P)$ 
should be open discs -- the \emph{faces} -- and the components of $S(P)\setminus V(P)$ 
should be open
segments -- the \emph{edges}. A standard spine is dual to a 1-vertex triangulation of $M$,
and this partially explains why $c(M)$ equals the minimal number of
tetrahedra in a triangulation when $M$ is \ptwoirred\ and distinct
from $S^3, \matRP^3, L_{3,1}$.

\paragraph{Sketch of the proof}
A closed non-orientable 3-manifold has a non-trivial first Stiefel-Whitney class $w_1\in
H^1(M;\matZ_2)$. A surface $\Sigma\subset M$ which is
Poincar\'e dual to $w_1$ is usually called a
\emph{Stiefel-Whitney surface}~\cite{HeLaNu}.
It has odd intersection with a loop $\gamma$ if and only
if $\gamma$ is orientation-reversing. 
It follows that
$M\setminus\Sigma$ is connected and orientable, \emph{i.e.}~$M = N \cup \calR(\Sigma)$ is
obtained by gluing a regular neighborhood $\calR(\Sigma)$ of $\Sigma$
to an orientable connected compact $N$ along their boundaries.

We can now list the main steps of the proof. Let $M$ be a
non-orientable \ptwoirred\ closed 3-manifold $M$ with $c(M)\leqslant 6$.
\begin{itemize}
\item[(1)] We prove that, without loss of generality, $\Sigma\subset M$ can be
assumed to lie in a minimal skeleton $P$ of $M$ so that $\calR(\Sigma)\cap P$
(whence $\calR(\Sigma)$) has some definite shape;

\item[(2)] using the shape of $\calR(\Sigma)\cap P$ we prove that $N$, with 
a suitable extra structure (a marking on $\partial N$) has a very low
(suitably defined) complexity;

\item[(3)] manifolds with marked boundary of low complexity are classified
in~\cite{MaPe}, so we list the possible shapes for $N$;

\item[(4)] we examine by hand how $\calR(\Sigma)$ and $N$ can be glued along 
$\partial\calR(\Sigma) = \partial N$, proving that precisely the four flat non-orientable 
manifolds and the torus bundle over $S^1$ of type {\rm Sol} with
monodromy trace $\matr 1110$ can arise;

\item[(5)] we exhibit some spines of manifolds of type $\matH^2\times S^1$, of type Sol,
and with non-trivial JSJ decomposition with $7$ vertices.

\end{itemize}

Our results on $\calR(\Sigma)\cap P$ are stated in the rest of this section and
proved in Section~\ref{proof_lemmas:sec}. The theory of complexity for manifolds
with marked boundary is reviewed in Section~\ref{teoria:section}, and is used
in Section~\ref{fine_dimo:sec}
to prove that $N$ has low complexity, and hence a definite shape.
The possible gluings
of $N$ and $\calR(\Sigma)$ are then analysed at the end of
Section~\ref{fine_dimo:sec}, 
to conclude the proof.

\paragraph{First part of the proof}
Let us start with a general result on Stiefel-Whitney surfaces.
\begin{prop} \label{prima:prop}
Let $\Sigma\subset M$ be a Stiefel-Whitney surface of a closed
non-orientable $M$.
The surfaces $\Sigma$ and $\partial \calR(\Sigma)$ are orientable. 
If $M$ is \ptwoirred\ then:
\begin{itemize}
\item $N = \closure{M\setminus\calR(\Sigma)}$ is \ptwoirred;

\item no component of $\Sigma$ or $\partial \calR(\Sigma)$ is a sphere;

\item if a component of $\Sigma$ or $\partial \calR(\Sigma)$ 
is a torus then it is incompressible.
\end{itemize}
\end{prop}
\begin{proof}
We first prove that $\Sigma$ is orientable. Suppose $\gamma\subset\Sigma$ is an
orientation-reversing loop (in $\Sigma$).
If $\gamma$ is
orientation-preserving in $M$ it can be perturbed to a loop
intersecting $\Sigma$ in one point, and if it is orientation-reversing
in $M$ it can be isotoped away from $\Sigma$:
both cases being in contrast to the definition of $\Sigma$.
Obviously, $\partial \calR(\Sigma)$ is orientable because $N$ is. 

Suppose now $M$ is \ptwoirred. Since $N$ is connected, each component of $\Sigma$ is
non-separating, thus it cannot be a sphere or a compressible torus. 
So no component of $\partial \calR(\Sigma)$ is a sphere. 
Suppose a component of $\partial \calR(\Sigma)$ is a compressible torus.
Then the corresponding component of $\calR(\Sigma)$ is the non-orientable
interval bundle $T\timtil I$ over the torus.
It follows quite easily that $M$ is a Dehn filling on $T\timtil I$, 
hence $S^2\timtil S^1$ or $\matP^2\times S^1$, 
a contradiction.

Let $S\subset N$ be a sphere. Then $S$ 
bounds in $M$ a ball, which cannot contain components of $\Sigma$ 
because they are non-separating. Hence the ball is contained in $N$, 
and the orientable $N$ is \ptwoirred.
\end{proof}

Let $P$ be a standard spine of a non-orientable $M$. The embedding $P\subset M$ induces
an isomorphism $H_2(P;\matZ_2)\cong
H_2(M;\matZ_2)$. Using cellular homology, a representative for a cycle in $H_2(P;\matZ_2)$ 
is a subpolyhedron consisting of some faces, an even number of them
incident to each edge of $P$. Such a subpolyhedron is a surface
near the edges it contains, and it is also a surface near the vertices (in fact, the link
of a vertex does not contain two disjoint circles).
Thus every
homology class is represented by a (unique) surface in $P$: in
particular there is a unique Stiefel-Whitney surface 
$\Sigma$ inside $P$.

Let us now suppose $M$ is \ptwoirred\ with $c(M)\leqslant 6$,
and $P$ is a minimal standard spine of $M$. 
The Stiefel-Whitney surface $\Sigma\subset P$ is not necessarily
connected but, since $P$ has at most 6 vertices,
it contains a few components of low genus.
Namely, we have the following result which will
be proved in Section~\ref{proof_lemmas:sec}.

\begin{lemma}\label{general_sigma:lem}
Let $M$ be \ptwoirred\ with $c(M)\leqslant 6$ and $P$ be a minimal standard spine of $M$.
The Stiefel-Whitney surface $\Sigma\subset P$ contains at most $2$ connected components.
Moreover $M$ has a minimal standard spine (which we denote again by $P$) 
with a Stiefel-Whitney surface (which we denote again by
$\Sigma$) having Euler characteristic equal to zero.
\end{lemma}

So $\Sigma$ consists of one or two tori. We fix a sufficiently small regular neighborhood
$\calR(\Sigma)$ of $\Sigma$ in $M$, such that the intersection of 
$\calR(\Sigma)$ and $P$ is a regular neighborhood $\calR_P(\Sigma)$ of $\Sigma$ in
$P$.
Using the fact that $P$ has 6 vertices at most, we will prove in
Section~\ref{proof_lemmas:sec} the following results. Recall that there
are two interval bundles on the torus up to homeomorphism, namely $T\times
I$ and $T\timtil I$.

\begin{figure}
\begin{center}
\mettifig{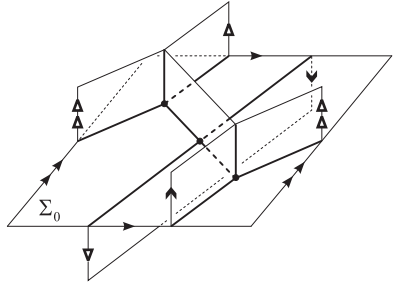} 
\caption{The regular neighborhood (in $P$) of a component $\Sigma_0$ of $\Sigma$, 
such that $\calR(\Sigma_0)=T\widetilde\times I$ (similar arrows must be identified).} 
\label{nhbd_T_non_orient:fig}
\end{center}
\end{figure}

\begin{lemma}\label{2_tori_sigma:lem}
Let $M$ be \ptwoirred\ with $c(M)\leqslant 6$ and $P$ be a minimal standard spine of $M$.
If $\Sigma\subset P$ consists of two tori, $\calR(\Sigma)$ consists of two copies 
of $T\timtil I$ and each of the components of $\calR_P(\Sigma)$
is as shown in Fig.~\ref{nhbd_T_non_orient:fig}.
\end{lemma}

\begin{figure}
\begin{center}
\mettifig{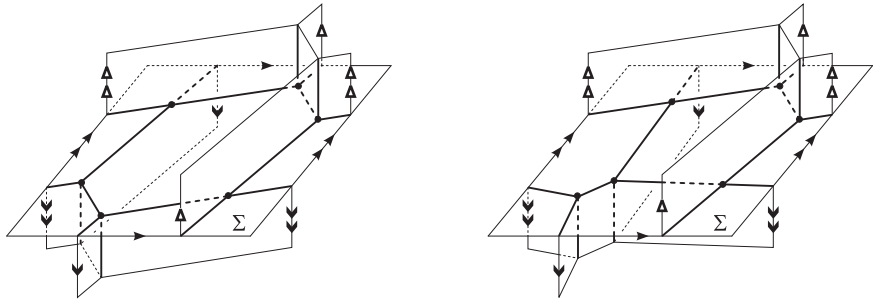} 
\caption{Two possibilities for the regular neighborhood (in $P$) of $\Sigma$, 
if $\Sigma$ consists of one torus and $\calR(\Sigma)=T\times I$
(similar arrows must be identified).}
\label{nhbd_T_orient:fig}
\end{center}
\end{figure}

\begin{lemma}\label{product_sigma:lem}
Let $M$ be \ptwoirred\ with $c(M)\leqslant 6$ and $P$ be a minimal standard spine of $M$.
If $\Sigma\subset P$ is one torus and $\calR(\Sigma)=T\times I$,
then $\calR_P(\Sigma)$ is one of the two polyhedra shown in Fig.~\ref{nhbd_T_orient:fig}.
\end{lemma}

\begin{lemma}\label{twisted_sigma:lem}
Let $M$ be \ptwoirred\ with $c(M)\leqslant 6$ and $P$ be a minimal standard spine of $M$.
If $\Sigma\subset P$ is one torus and
$\calR(\Sigma)=T\timtil I$, then $M$ has a
minimal standard spine
(which we denote again by $P$) with a Stiefel-Whitney
surface $\Sigma$ such that $\calR(\Sigma) = T\timtil I$ and
$\calR_P(\Sigma)$ is as shown in Fig.~\ref{nhbd_T_non_orient:fig}.
\end{lemma}

We now know that $\calR_P(\Sigma)$ has 3 possible shapes.
In order to complete our classification, we need to know the possible
shapes of the rest of $P$, namely the polyhedron
$Q=\closure{P\setminus\calR_P(\Sigma)}$. Moreover, we know that the two
polyhedra are glued along a very special graph contained in
$\partial\calR(\Sigma)$:
it consists of either one or two $\theta$-graphs.
Here, a $\theta$-graph is a trivalent 
graph $\theta$ contained in a torus $T$ such that $T\setminus\theta$
is an open disc.
Decompositions of minimal spines (and manifolds) along $\theta$-graphs (and tori)
have been studied in~\cite{MaPe, MaPe:nonori}. The basic
result is a decomposition theorem for \ptwoirred\ manifolds. Since we will use
it on $N$, which is orientable, we describe in Section~\ref{teoria:section}
the orientable version of the theory. 
(We only note here that non-orientable 
manifolds could be cut along Klein bottles
also, and the graph \mettifig{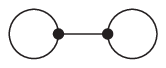, width=1 cm} 
should be taken into account in this case.)
We will then conclude the proof
of Theorem~\ref{main:teo} in Section~\ref{fine_dimo:sec}.

\section{Manifolds with marked boundary}\label{teoria:section}

\begin{figure}
\begin{center}
\mettifig{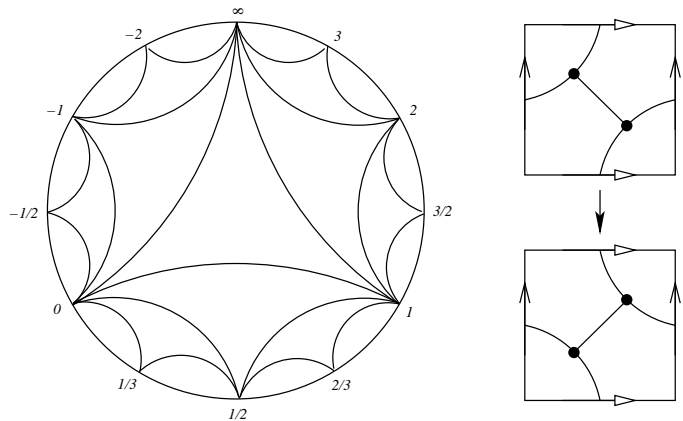} 
\caption{A tesselation of the Poincar\'e disc into ideal triangles (left) and a flip (right).}
\label{tesselation:fig}
\end{center}
\end{figure}

\paragraph{$\theta$-graphs in the torus}
A {\em $\theta$-graph} in the torus $T$ is a trivalent graph $\theta\subset T$  such that
$T\setminus\theta$ is an open disc. The embedding of $\theta$ in $T$
is unique up to homeomorphism of $T$, but \emph{not} up to
isotopy. There is a nice description, taken from~\cite{FloHa},
of all $\theta$-graphs (up to
isotopy) which we now describe. After fixing a basis $(a,b)$ for
$H_1(T;\matZ)$, every \emph{slope} on $T$ (\emph{i.e.}~isotopy class of simple closed
essential curves) is determined by its unsigned homology class $\pm(pa+qb)$,
thus by the number $p/q\in\matQ\cup\{\infty\}$. Consider
$\matQ\cup\{\infty\}$ sitting inside $\matR\cup\{\infty\}$, the
boundary of the upper half-plane of $\matC$, with its standard
hyperbolic metric. For each pair $p/q, r/s$ of slopes having algebraic
intersection $\pm 1$ (\emph{i.e.}~such that $ps-qr=\pm 1$) draw a
geodesic connecting $p/q$ and $r/s$. The result is a tesselation of
the half-plane into ideal triangles, shown in Fig.~\ref{tesselation:fig}-left (in the disc
model).

It is easily seen that a $\theta$-graph is determined by the three
slopes it contains, and that such slopes have pairwise
intersection 1. Thus, a $\theta$-curve corresponds to a triangle of
the ideal tesselation, \emph{i.e.}~to a vertex of the dual trivalent tree.
Moreover, two $\theta$-graphs are connected by a segment in this tree
when they share two slopes, \emph{i.e.}~when we can pass from one
$\theta$-graph to the other via a \emph{flip}, shown in Fig.~\ref{tesselation:fig}-right.

\paragraph{Manifolds with marked boundary}
Let $M$ be a connected compact 3-manifold with (possibly empty)
boundary consisting of tori. By associating to each torus component of
$\partial M$ a $\theta$-graph, we get a
{\em manifold with marked boundary}. As we have seen, the same manifold can
be marked in infinitely many distinct ways.

Now we describe two fundamental operations on the set of manifolds with marked boundary.
The first one is binary:
if $M$ and $M'$ are two such objects, take two tori $T\subset\partial
M, T'\subset\partial M'$ marked with $\theta\subset T, \theta'\subset
T'$ and a homeomorphism $\psi:T\stackrel{\sim}{\rightarrow} T'$ such that $\psi(\theta)=\theta'$.
By gluing $M$ and $M'$ along $\psi$ we get a new 3-manifold with marked
boundary. We call this operation an \emph{assembling}.
Note that, although there are infinitely many non-isotopic maps between
two tori, only finitely many of them send one marking to the other, so
there is a finite number of non-equivalent assemblings of $M$ and $M'$.

We describe the second operation.
Let $M$ be a manifold with marked boundary, and $T,T'$ be two
distinct boundary components of it, marked with $\theta\subset T$ and $\theta'\subset
T'$.
Let $\psi : T \stackrel{\sim}{\rightarrow} T'$ be a homeomorphism such
that $\psi(\theta)$ equals either $\theta'$ or a $\theta$-graph obtained
from $\theta'$ via a flip. The manifold obtained identifying $T$ and $T'$
via $\psi$ is a new manifold with marked boundary. (There is a
technical reason for not asking only that $\psi(\theta)=\theta'$,
which will be clear later.)
We call this operation a \emph{self-assembling}.
Again, there is only a finite number of non-equivalent self-assemblings.

\paragraph{Spines and skeleta} 
The notion of spine extends to the class of manifolds with marked
boundary. 
A sub-polyhedron $P$ of a 3-manifold with marked boundary $M$ is
called a {\em skeleton} 
of $M$ if $M\setminus(P \cup \partial M)$ is
an open ball and $P \cap \partial M$ is a graph contained in the
marking of $\partial M$.
We have not used the word ``spine'' because maybe $P$ is not a spine of $M$
in the usual sense when $\partial M\neq\emptyset$ 
-- \emph{i.e.}~$M$ does not retract onto $P$. 
On the other side note that, if $M$ is closed, 
a skeleton of $M$ is just a spine of $M$.
Recall that a polyhedron is \emph{simple} when the link of every point
is contained in the 1-skeleton of a tetrahedron $K$. 
It is easy to prove that each 3-manifold with marked boundary has a
simple skeleton.

\paragraph{Complexity}
The {\em complexity} of a 3-manifold with marked boundary $M$ is of
course defined as the minimal number of vertices of a simple skeleton of 
$M$.
It depends on the topology of $M$ and on the marking. In particular, if $T=\partial
M$ is one torus then every (isotopy class of a) $\theta$-graph on $T$
gives a distinct complexity for $M$.
Three properties extend from the closed case to the case with marked boundary:
complexity is still additive under connected sums, it is finite-to-one on
orientable irreducible manifolds with marked boundary,
and if $M$ is orientable irreducible with $c(M)>0$, then it has a
minimal standard skeleton~\cite{MaPe}.
A skeleton $P\subset M$ is called standard when $P\cup\partial M$ is.

\begin{figure}
\begin{center}
\mettifig{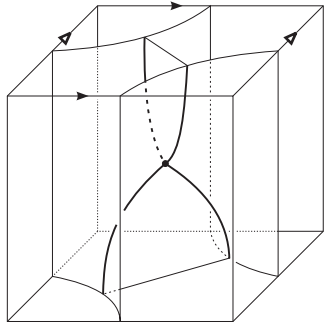} 
\caption{A skeleton for $B_3$.} \label{skeleton_B3:fig}
\end{center}
\end{figure}

\paragraph{Examples}
Let $T$ be the torus. Consider $M=T\times I$, the boundary being
marked with a $\theta_0\subset
T\times 0$ and a $\theta_1\subset T\times 1$. If $\theta_0$ and
$\theta_1$ are isotopic, the resulting 
manifold with marked boundary is called $B_0$. If $\theta_0$ 
and $\theta_1$ are related by a flip, we call the resulting manifold
with marked boundary $B_3$.
A skeleton for $B_0$ is $\theta_0 \times [0,1]$, 
while a skeleton for $B_3$ is shown in
Fig.~\ref{skeleton_B3:fig}. The skeleton of $B_0$ has no vertices, so $c(B_0)=0$. The
skeleton of $B_3$ has 1 vertex, and it can be
shown~\cite{MaPe} that there is no skeleton for $B_3$ without vertices, so $c(B_3)=1$.

\begin{figure}
\begin{center}
\mettifig{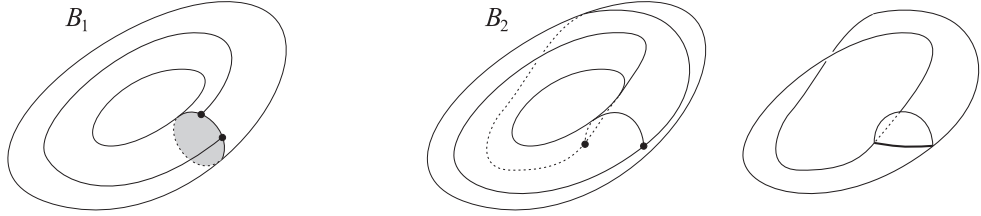, width = 13 cm} 
\caption{The manifolds with marked boundary $B_1$ (left) and $B_2$ (centre), and a skeleton
for $B_2$ (right).}
\label{B1_B2:fig}
\end{center}
\end{figure}

Two distinct marked solid tori are shown in Fig.~\ref{B1_B2:fig} (left and centre) and denoted
by $B_1$ and $B_2$. A skeleton for $B_1$ is a meridinal disc with boundary contained in
the $\theta$-graph. A skeleton for $B_2$ is shown in Fig.~\ref{B1_B2:fig}-right.
Since they have no vertices, we have $c(B_1)=c(B_2)=0$.

\begin{figure}
\begin{center}
\mettifig{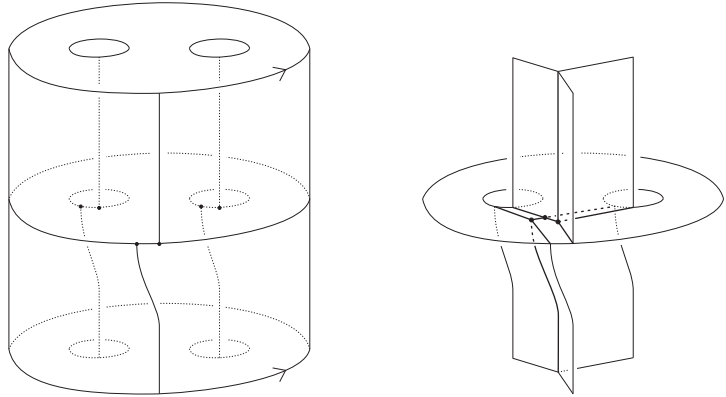}
\caption{The manifold with marked boundary $B_4$ (left) and a skeleton for it (right).} 
\label{B4:fig}
\end{center}
\end{figure}

The first irreducible orientable manifold with more than two marked
boundary components has $c=3$.
Let $D_2$ be a disc with two holes. Set $M=D_2\times S^1$. For each
torus $T$ in $\partial M$, a basis $(a,b)$ for $H_1(T;\matZ)$ is given by
taking $a$ to be $\partial D_2\times \{{\rm pt}\}$
(with orientation induced from that of
$D_2$) and $b$ to be $\{{\rm pt}\}\times S^1$, oriented as $S^1$. 
With respect to this basis, on each boundary component 
a triple of slopes $\{i,\infty,i+1\}$ defines a $\theta$-graph
$\theta_i$ for any integer $i$ (note that $\theta_{-1}$ and $\theta_0$ are the
$\theta$-graphs containing $0$ and $\infty$). 
Now let $B_4$ be $M$ with markings $\theta_0, \theta_0$
and $\theta_{-1}$, see Fig.~\ref{B4:fig}-left.
It has a skeleton with 3 vertices, shown in Fig.~\ref{B4:fig}-right. 
It can be proved~\cite{MaPe} that $B_4$ has no skeleton with less vertices, so $c(B_4)=3$,
and that a distinct choice for the markings -- for instance the same $\theta_0$ on all
boundary components -- would give $c>3$.

\paragraph{Assemblings and skeleta}
Let $M, M'$ be two manifolds with marked boundary, and $P, P'$ be two
corresponding standard skeleta. An assembling of $M$
and $M'$ is given by a map $\psi$ that matches the $\theta$-graphs, so
$P\cup_\psi P'$ is a simple polyhedron inside $M\cup_\psi M'$. 
Moreover, it is not difficult to see that $P\cup_\psi P'$ is a
skeleton of the new manifold with marked boundary $M\cup_\psi M'$.
(This is true because the complement of a $\theta$-graph is a single
disc, so the complement of $P\cup_\psi P'$ consists of two balls glued
together along a single disc, hence another ball.)

If $P,P'$ have $n, n'$ vertices, then $P\cup_\psi P'$ has $n+n'$
vertices. Suppose $P$ and $P'$ are minimal skeleta of $M$ and $M'$,
\emph{i.e.}~$n$ and $n'$ are the complexities of $M$ and $M'$. It is not
true in general that $P\cup_\psi P'$ is minimal.
Since $M\cup_\psi M'$ has a skeleton with $n+n'$ vertices, its
complexity is at most $n+n'$, and it equals $n+n'$ precisely when $P\cup_\psi
P'$ is minimal.
We will be interested in the case when $P\cup_\psi P'$ is
minimal: in other words, complexity is sub-additive
under assemblings, 
and we will be interested in the case when it is additive.

An analogous construction works for self-assemblings. Let $M'$ be
obtained self-assembling $M$, along a map $\psi:T\stackrel{\sim}{\to}T'$ such that
$\psi(\theta)$ either equals $\theta'$ or is obtained from $\theta'$
via a flip. In any case,
it is possible to isotope $\psi$ to $\psi'$ so that $\psi'(\theta)$ and $\theta'$
intersect each other transversely in 2 points, and to use the map
$\psi'$ to construct $M'$.
Let $P$ be a standard
skeleton for $M$.
Take $P'=P\cup T$ inside $M'$:
again, $P'$ is a skeleton for $M'$. The polyhedron $P'$ is the result
of adding one of the two polyhedra shown in
Fig.~\ref{nhbd_T_orient:fig} to $P$.
(Note that a construction
analogous to the one made for assemblings does not work:  if
$\psi(\theta)=\theta'$, then $P$ alone inside $M'$ is not a
skeleton of $M'$, because its complement is a solid torus: this is why
it is necessary to add $T$. Moreover $P\cup T$ is a skeleton but is
not standard, so we need to isotope $\psi$ to recover
standardness.) 
This operation creates 6 new vertices: 
if $P$ has $n$ vertices, then $P'$ has $n+6$ vertices. So the
complexity of $M'$ is at most the complexity of $M$ plus 6.

\paragraph{Bricks}
The theory ends with a decomposition theorem. An assembling is
\emph{sharp} if the complexity is additive and both manifolds with
marked boundary are irreducible and distinct from $B_0$, and a self-assembling is
\emph{sharp} if the complexity of the new manifold is the complexity
of the old one plus 6.
An irreducible orientable manifold with marked boundary is a {\em brick} if it is not the
result of a sharp assembling or self-assembling of other irreducible
manifolds with marked boundary. The proof of the following
result is clear: if an irreducible manifold with marked boundary is not a brick, then it can be
de-assembled. Then we repeat the analysis on each new piece. 
Since the sum of the complexities of all pieces does not increase
(and since the only possible pieces with complexity 0 are known to be
$B_1$ and $B_2$),
this iteration must stop after finite time.
\begin{prop} \label{decomposition:prop}
Every irreducible orientable manifold with marked boundary can
be obtained from some 
bricks via a combination of sharp assemblings and sharp self-assemblings.
\end{prop}
This result can be restated at the level of skeleta: 
every orientable manifold with marked boundary has a minimal skeleton which
splits along $\theta$-graphs into minimal skeleta of bricks. Here, bricks
are defined to be orientable. (Non-orientable bricks are analogously 
defined in~\cite{MaPe:nonori}, but
we do not need them here.)

It is proved in~\cite{MaPe} that the only bricks with boundary having
complexity at most 3 are the $B_0,\ldots, B_4$ introduced above. Using
a computer, all bricks $B_0,\ldots,B_{10}$ having complexity up to 9
and with non-empty boundary
have been classified. Let $P_i$ be a minimal skeleton of $B_i$: 
Proposition~\ref{decomposition:prop} implies that every orientable manifold
having complexity at most $9$ has a minimal spine which splits along
$\theta$-graphs into copies of $P_1, \ldots, P_{10}$.
Bricks $B_5,\ldots, B_{10}$ have complexity $8$
or $9$, moreover they are all hyperbolic except $B_5$.

\paragraph{Assembling small bricks}
Let $M$ be a manifold with marked boundary. Let us examine the effect
of assembling $M$ with some $B_i$ along a torus $T\subset\partial M$,
marked with a $\theta\subset T$. 
Choose a basis for $H_1(T;\matZ)$ so that $\theta$ corresponds to the triple 
$\{0,1,\infty\}$, see Fig.~\ref{tesselation:fig}-left.
If $i=0$, the
assembling leaves $M$ unaffected. If $i=1$, a 
Dehn filling is performed on $M$, killing one of the three
slopes $0,1,\infty$. 
If $i=2$, a Dehn filling is performed on $M$, killing one of the slopes $2, 1/2, -1$.
If $i=3$, the graph $\theta$ is changed by a flip. It follows that by
assembling $M$ with some copies of $B_1, B_2,$ and $B_3$ we can
arbitrarily change some markings or do arbitrary Dehn fillings on $M$.

We can use Proposition~\ref{decomposition:prop} and the known list
of bricks to classify manifolds with non-empty marked boundary of low
complexity. Every such manifold is obtained via \emph{sharp}
assemblings and self-assemblings from the known bricks. For instance, if
a marked $M$ has complexity at most 2, no self-assembling is involved
since it adds 6 to the complexity, and only assemblings of $B_0, B_1,
B_2,$ and $B_3$ are involved. Therefore $M$ is a (marked) solid torus,
or a (marked) product $T\times I$. We are here interested in the first case where
$M$ has one boundary component and is not a (marked) solid torus.
Let $(D_2\times S^1)_{2,2,\theta_{-1}}$ be the Seifert manifold with
base space a disc and two fibers of type $(2,1)$, marked with 
$\theta_{-1}$ in the boundary. (Recall that 
$\theta_{-1}$ is the $\theta$-graph containing
the slopes $\infty,-1,$ and $0$, where coordinates are taken with respect to
the obvious basis of $H_1(D_2\times S^1;\matZ)$.)

\begin{prop} \label{only:3:prop}
Every irreducible manifold with a single marked boundary component 
having $c\leqslant 2$ is a marked solid torus.
Every such manifold having $c=3$        
is a marked solid torus or
$(D_2\times S^1)_{2,2,\theta_{-1}}$.
\end{prop}
\begin{proof}
Suppose $M$ is not a marked solid torus, with $c\leqslant 3$.
It decomposes into copies of $B_1, B_2, B_3,$ and $B_4$,
and at least one $B_4$ must be
present. Moreover, since $c(B_4)=3$, the other
bricks in the assembling have complexity 0, so they must be
$B_1$'s and $B_2$'s. Despite the
apparent lack of symmetry of the markings, for each pair of boundary
components there is an automorphism of $B_4$ interchanging them (and
their markings), so it is not important to which boundary components
the $B_1$'s and $B_2$'s are assembled. 
Suppose then the assemblings are performed on
the first two components. It follows from the discussion above that
we can realize Dehn fillings on slopes
$\infty, 0, 1$ with $B_1$ and on slopes $-1,1/2,2$ with $B_2$. The
only such filling that creates a singular fiber is $2$, whence the result.
\end{proof}

\paragraph{A manifold with non-trivial JSJ decomposition containing hyperbolic pieces}
It is now easy to use the known bricks to build manifolds. Using
$B_1, B_2, B_3,$ and $B_4$ any graph manifold can be built. The brick
$B_6$ is the first hyperbolic brick, having $c=8$ (whereas
$B_7,\ldots, B_{10}$ have $c=9$, see~\cite{MaPe}). It is the
figure-eight knot sister, denoted by $M2^1_2$ in~\cite{CaHiWe}, marked
with the most natural $\theta$-graph: it is the $\theta$-graph
containing the 3 shortest slopes in the cusp, or equivalently the
unique $\theta$-graph fixed by any isometry of $M2^1_2$.
Note that any other marked hyperbolic manifold has $c>8$: the manifold
$M2^1_2$ is then in some sense `smaller' (or `simpler') than the figure-8
knot complement $M2^1_1$, although they have the same volume -- note
that the smallest known closed hyperbolic manifold is obtained via Dehn
filling from $M2^1_2$ but not from $M2^1_1$.

If we assemble $B_6$ with $B_1$ or $B_2$, we always get a
non-hyperbolic manifold: in order to get a hyperbolic one, we must use
a $B_3$ and a $B_2$, which is coherent with the fact that the first
closed hyperbolic manifolds have $c=9=8+1=c(B_6)+c(B_3)$.
It is easier to construct a closed manifold whose JSJ decomposition is non-trivial
and contains a hyperbolic piece: simply take any assembling of $B_6$
and $(D_2\times S^1)_{2,2,\theta_{-1}}$. The complexities of the pieces
are $8$ and $3$, so we get a manifold with $c\leqslant 11$, but we
cannot be sure that equality holds -- in other words, by gluing
minimal spines of $B_6$ and $(D_2\times S^1)_{2,2,\theta_{-1}}$
we get a spine of the closed
manifold, which is possibly not minimal.
Nevertheless, we know from~\cite{MaPe} that every brick with $c\leqslant 9$ is atoroidal. 
If this were true for any $c$, every sharp decomposition of
a closed irreducible manifold into bricks would be a refinement of its
JSJ decomposition. In
other words, there would be a minimal spine of the closed manifold which decomposes
into minimal skeleta of the pieces of the JSJ decomposition (which might
further decompose), with appropriate markings.
Therefore, the complexity of a closed manifold would be the sum of the
complexities of the (appropriately marked) pieces of its JSJ decomposition: in
particular, a hyperbolic piece would give a contribution $\geqslant 8$, 
and a Seifert one a contribution $\geqslant 3$, giving $8+3=11$ at least.

\section{End of the proof} \label{fine_dimo:sec}

At the end of Section~\ref{complexity:section}, we have listed the
possible shapes for the regular neighborhood $\calR_P(\Sigma)$ of the Stiefel-Whitney
surface $\Sigma$ in $P$. In all cases, the polyhedron $P$ can be cut
into a (possibly disconnected) $\calR_P(\Sigma)$ and a connected 
$Q=\closure{P\setminus\calR_P(\Sigma)}$. The two subpolyhedra are glued
along $\theta$-graphs. At the level of manifolds, $\calR_P(\Sigma)$ is
contained in $\calR(\Sigma)$ and $Q$ is contained in $N$, which is
orientable and irreducible.
The original \ptwoirred\ $M$ is decomposed along one or two tori into
$\calR(\Sigma)$ and $N$.
Both $\calR(\Sigma)$ and $N$ are equipped with a marking
on each boundary component,
given by the $\theta$-graphs separating the polyhedra.
It is easy to check that $N\setminus(\partial N\cup Q)$ is an open
ball in any case, so $Q$ is a skeleton of $N$. Concerning the 3
possible shapes for $\calR_P(\Sigma)$, two of them are skeleta of the
corresponding $\calR(\Sigma)$, and the other one is not. We now study
this in detail.

\paragraph{If $\Sigma$ consists of two tori}
Each component of $\calR_P(\Sigma)$ is a skeleton (with 3 vertices)
of the corresponding marked
$T\timtil I$. Therefore $M$ is obtained assembling a copy of
the marked $T\timtil I$ on each boundary component of $N$. 
Since $\calR_P(\Sigma)$ has 6 vertices and $P$ has $6$ vertices at
most, there is no vertex in $Q$,
so $N$ has complexity zero (and it is orientable), hence $N=B_0$ is the trivial brick. 
Thus $M$ is obtained assembling two copies of the marked
$T\timtil I$.

We now prove that the result of this assembling must be a flat manifold.
Note that $T\timtil I$ has two distinct
fibrations: the first one is the product $S\times S^1$, where $S$ is the
M\"obius strip. The second one is a Seifert fibration over the
orbifold whose underlying topological space is an annulus, with one
mirror circle (so the orbifold has only one true boundary component, see~\cite{Sco}).
A basis $(a,b)$ for $H_1(\partial(S\times S^1);\matZ)$ is
given by taking $a=\partial S\times\{{\rm pt}\}$ and $b=\{{\rm pt}\}\times S^1$, 
with some orientations. With respect to this basis,
slopes are numbers in $\matQ\cup\{\infty\}$, and
$0$ is a fiber of the first fibration, while $\infty$ is a fiber of
the second fibration. The $\theta$-graph in the
boundary is the one containing the slopes $\infty,0,1$. An assembling
of two copies of $S\times S^1$ is given by a map $\psi$ which
matches the markings, \emph{i.e.}~sends the set of slopes $\{\infty,0,1\}$ of
the first one to the set $\{\infty,0,1\}$ of the other one.

If $\psi(0)=0$, we get a fibration over the Klein bottle. If
$\psi(0)=\infty$ or $\psi(\infty)=0$, 
we get a fibration over a M\"obius strip with one mirror
circle. If $\psi(\infty)=\infty$, we get a fibration over an annulus with two
mirror circles. In all cases the base orbifold has
$\chi^{\rm orb}=0$, so the manifold is flat.

\paragraph{If $\Sigma$ is one torus and $\calR(\Sigma)=T\widetilde\times I$}
The polyhedron $\calR_P(\Sigma)$
is a skeleton of the marked $T\timtil I$. Therefore $M$ is
obtained by assembling $N$ with $T\timtil I$.
Since $\calR_P(\Sigma)$ has 3 vertices, there are three vertices at most in
$Q$.
Proposition~\ref{prima:prop} implies that $N$ is not a (marked) solid
torus, hence 
$N = (D_2\times S^1)_{2,2,\theta_{-1}}$ by Proposition~\ref{only:3:prop}.

As above, we prove that the result of this assembling must be a flat
manifold. Note first that $(D_2\times S^1)_{2,2,\theta_{-1}}$ fibers
over a disc with two singular fibers of type $(2,1)$, or as a twisted
product $S\timtil S^1$ over the M\"obius strip $S$.
The $\theta$-curve $\theta_{-1}$ contains the slopes $\infty,-1,0$,
and $0$ is a fiber of the first fibration, while $-1$ is a fiber of
the second fibration. Now, an assembling is given by a map $\psi$ that
sends the triple of slopes $\{\infty,-1,0\}$ of $N$ to the triple of
slopes $\{\infty,0,1\}$ of $T\timtil I$.
If $\psi(0)=0$, we get a
fibration over $\matRP^2$ with two singular fibers of type $(2,1)$. If
$\psi(0)=1$ we get a fibration over a disc with two singular fibers of
type $(2,1)$ and a mirror circle. If $\psi(-1)=0$ we get a fibration
over a Klein bottle, and if $\psi(-1)=1$ we get a fibration over a
M\"obius strip with one mirror circle. In all cases the 
base orbifold has $\chi^{\rm orb}=0$, so the manifold is flat.

\paragraph{If $\Sigma$ is one torus and $\calR(\Sigma)=T\times I$}
The polyhedron $\calR_P(\Sigma)$ is \emph{not} a skeleton of the
marked $T\times I$, since $(T\times I)\setminus(\partial(T\times
I)\cup\calR_P(\Sigma))$ consists of two balls instead of one. 
The polyhedron $\calR_P(\Sigma)$ is one involved in self-assemblings,
so our $M$ is the result of a self-assembling of
$N$, which has two boundary components and complexity 0 (because 6
vertices are in $\calR_P(\Sigma)$). 
Therefore $M$ is a self-assembling
of $B_0$, \emph{i.e.}~it is the mapping torus of a map $\psi:T\to T$ that sends
a $\theta$-graph $\theta\subset T$ to a $\theta$-graph 
$\psi(\theta)$ sharing at
least two slopes with $\theta$. Let $a,b\in H_1(T;\matZ)$ represent these
two slopes.
With respect to the basis $(a,b)$, we have
$\psi(0),\psi(\infty)\in\{0,1,\infty\}$, therefore $\psi$ is read
as a matrix $A$, with trace between $-1$ and $1$ (recall that $M$ is
non-orientable), which is either periodic or hyperbolic.
In the former case $M$ is flat, while in the latter one $M$ is 
a torus bundle (of type Sol) with monodromy $A$
(see~\cite{Sco}).
Obviously, these self-assemblings are sharp because $P$ is minimal.
Now, it is easy to prove that every non-periodic matrix $A\in{\rm GL}_2(\matZ)$ with
$\det A=-1$ and $|{\rm tr}\, A|\leqslant 1$ is conjugated either to $\matr
1110$ or to $\matr 1110^{{\tiny -1}}=\matr 011{-1}$ (so ${\rm tr}\, A=\pm
1$), hence there is only one such manifold of type Sol.

\paragraph{Conclusion}
We have proved that every closed non-orientable \ptwoirred\ manifold
with $c\leqslant 6$ either is flat or it is the
torus bundle (of type {\rm Sol}) with monodromy $\matr 1110$, and that
it has complexity 6.
Moreover, each of these 5 manifolds occurs.
In fact, each of the 4 flat manifolds fibres in a few distinct ways
over 1-~or 2-dimensional orbifolds, and it follows
from~\cite{Sco} that all 4 can be realized with some
of the fibrations described above.
Moreover, the Sol manifold has been constructed by self-assembling
$B_0$ sharply.

\paragraph{Examples of non-orientable manifolds with complexity 7}
Using $T\timtil I$, $B_4$, $B_3$, and two $B_2$'s, it is now 
easy to construct closed manifolds of complexity $7$. We have seen that a 
Dehn filling killing a slope in $\{\infty, 0, 1, -1,1/2,2\}$ on the first
component of $\partial B_4$ can be realized assembling $B_1$ or
$B_2$. 
We can kill
the slope $3$ as follows: we first assemble $B_3$, so that $\theta_0$ 
is replaced by $\theta_1$ (the $\theta$-graph corresponding to $\{1,2,\infty\}$), 
and then we assemble $B_2$. Therefore the manifold with marked boundary
$(D_2\times S^1)_{3,2,\theta_{-1}}$, obtained from $B_4$ by filling the first two boundary
components along the slopes $3$ and $2$, can be realized with a $B_4$, a $B_3$, and
two $B_2$'s, thus it has complexity at most $4$. Now we can assemble it with
the marked $T\timtil I$ considered above, along a map $\psi$.
The manifold $(D_2\times S^1)_{3,2,\theta_{-1}}$ has one fibration only, with
fiber $0$, whereas $T\timtil I$ has two, with fibers $0$ and $\infty$.
If $\psi(0)=0$, we get a Seifert fibration over $\matR\matP^2$ with
two critical fibres with Seifert invariants $(3,1)$ and $(2,1)$.
If $\psi(0)=\infty$, we get a Seifert fibration over the disc with one
reflector circle and two critical fibres with Seifert invariants
$(3,1)$ and $(2,1)$.
In both cases we get $\chi^{\rm orb}<0$, hence a manifold of type
$\matH^2\times\matR$.
Moreover, the complexity is $4+3=7$ at most, hence it is $7$.

A manifold of type Sol with $c=7$ can be easily constructed as above, with a 
self-assembling of $B_3$ realizing the monodromy $\matr 2110$ by sending
$\{\infty,0,1\}$ to $\{2,\infty,3\}$.

\section{Proofs of the lemmas}\label{proof_lemmas:sec}

We conclude with the proofs of the four lemmas of Section~\ref{complexity:section}.
First, we state (and prove) some easy properties of a minimal standard
spine of a non-orientable manifold with arbitrary number of vertices.
Then, we prove the lemmas.

The following criteria for non-minimality are proved in~\cite{MaPe,
MaPe:nonori}. Let $P$ be a standard spine of a closed \ptwoirred\ manifold. Then:
\begin{enumerate}

\item If a face of $P$ is embedded and incident to 3
or fewer vertices, $P$ is not minimal.
\item If a loop, embedded in $P$, intersects 
transversely the
singularity of $P$ in 1 point and bounds a disc in the 
complement of $P$, then $P$ is not minimal.
\end{enumerate}

\begin{figure}
\begin{center}
\mettifig{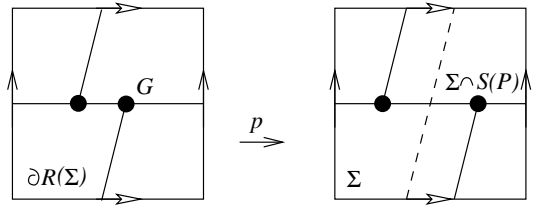} 
\caption{An example of map $p$.}
\label{example:fig}
\end{center}
\end{figure}

Throughout this section we suppose $P$ to be a minimal standard spine of a non-orientable
manifold $M$.
In the first part of this section we do not ask that $c(M)\leqslant 6$, so we allow $P$ to have an arbitrary number of vertices.
After, when we will prove the four lemmas, we will come back to the case when $P$ has at most $6$ vertices.
As above we call $\Sigma$ the Stiefel-Whitney surface of $M$ contained in $P$.
We fix a small regular neighborhood
$\calR(\Sigma)$ of $\Sigma$ in $M$, such that the intersection of 
$\calR(\Sigma)$ and $P$ is a regular neighborhood of $\Sigma$ in
$P$. We denote by $p:\partial \calR(\Sigma)\to\Sigma$ the projection. 
Then $(\calR(\Sigma)\setminus\Sigma)\cap P = G\times[0,1)$, where
$G = P\cap \partial\calR(\Sigma)$ is a trivalent graph.
The graph $p(G)$ has vertices with valence 3 and 4, and it 
is the intersection of $\Sigma$ and the singular set $S(P)$ of
$P$. The map $p:G\to \Sigma$
is a transverse immersion, \emph{i.e.}~it
is injective except in some pairs of points of $G$, that have the same
image, creating the 4-valent vertices of $p(G)$.
See an important example in 
Fig.~\ref{example:fig} with $\Sigma$ a torus and $\calR(\Sigma)=T\timtil I$.
The graphs $G$ and $p(G)$ fulfill some requirements, due to the
minimality of $P$.

\begin{lemma} \label{no:disc:lem}
No component of $G$ is contained in a disc of $\partial\calR(\Sigma)$.
\end{lemma}
\begin{proof}
Suppose a component $G_0$ is contained in a disc.
If $p$ is injective on $G_0$, then $p(G_0)$ is connected and contained
in a disc of $\Sigma$, but $\Sigma\setminus p(G)$ consists of discs
(because $P$ is standard), a contradiction.
If $p$ is not injective on $G_0$, then $p(G)$ intersects some edges of
$p(G\setminus G_0)$. But we can shrink and isotope $G_0$, and
consequently $G_0\times
[0,1)$ and $p(G_0)$, so that $p(G_0)$ does not intersect any edge of
$p(G\setminus G_0)$. The result is another spine of the same manifold,
but with fewer vertices: a contradiction.
\end{proof}

Since $P$ is standard, $\Sigma\setminus p(G)$ consists of
discs. Concerning $\partial\calR(\Sigma)$, we can only prove that
$\partial\calR(\Sigma)\setminus G$ can be embedded in the 
2-sphere and that it consists of discs inside the torus
components of $\partial\calR(\Sigma)$.

\begin{lemma} \label{sphere_minus_discs:lem}
The set $\partial\calR(\Sigma)\setminus G$ can be embedded in the
$2$-sphere.
\end{lemma}
\begin{proof}
The set $\partial\calR(\Sigma)\setminus G$ can be seen as a subset of the regular 
neighborhood $\calR(P)$ of $P$ in $M$ which is a sphere (because $M
\setminus P$ is a ball).
\end{proof}

\begin{lemma}\label{discs:lem}
Let $T$ be a torus component of $\partial\calR(\Sigma)$. 
Then $T\setminus G$ consists of discs.
\end{lemma}
\begin{proof}
Suppose a component $C$ of $T\setminus G$ contains a loop
$\alpha$ essential in $C$. Then $\alpha$ is essential in the whole of $\partial\calR(\Sigma)$
by Lemma~\ref{no:disc:lem}.
The loop $\alpha$ is unknotted in the ball $M\setminus P$, thus it
bounds a disc. Therefore $T$ is compressible in $N$, in contrast to 
Proposition~\ref{prima:prop}.
\end{proof}

Since $\Sigma\setminus p(G)$ consists of discs, connected components of
$p(G)$ correspond to connected components of $\Sigma$.

\begin{lemma}\label{one:vertex:lem}
Every connected component of $p(G)$ contains at least one $4$-valent vertex.
\end{lemma}
\begin{proof}
Let $\Sigma_0$ be a connected component of $\Sigma$.
The graph $p(G)\cap \Sigma_0$ is a connected 
component of $p(G)$ and $\Sigma_0\setminus (p(G)\cap \Sigma_0)$ is made of discs.
These two facts easily imply that if $p(G)\cap \Sigma_0$ contains only 3-valent
vertices then $P$ lays on a well-defined side of $\Sigma_0$, so we can choose
a transverse orientation for $\Sigma_0$.
Hence, $\partial\calR(\Sigma)\setminus G$ contains a surface
homeomorphic to $\Sigma_0$, contradicting
Lemma~\ref{sphere_minus_discs:lem}.
\end{proof}

\begin{lemma}\label{existence_3_valent_vertices:lem}
If $\Sigma$ is not connected, then every component of $p(G)$
contains at least one $3$-valent vertex.
\end{lemma}
\begin{proof}
Suppose a component of $p(G)$ contained in a component $\Sigma_0$ of
$\Sigma$ contains only 4-valent vertices.
Then $p(G)\cap\Sigma_0$ is a connected component of $S(P)$.
Since $S(P)$ is connected, then $p(G)\cap\Sigma_0 = S(P)$.
Obviously, each component of $\Sigma$ different from $\Sigma_0$ also
contains some singular points of $S(P)$.
A contradiction.
\end{proof}

\begin{lemma}\label{3_valent_two_discs:lem}
If a connected component of $p(G)$ (corresponding to a connected 
component $\Sigma_0$ of $\Sigma$) contains a $3$-valent vertex, 
then $\Sigma_0\setminus p(G)$ is
made of at least two discs.
\end{lemma}
\begin{proof}
We prove that the 3 germs of discs incident to a 3-valent vertex $v$
cannot belong to the same disc.
Suppose by contradiction that they do, and call $D$ this disc.
Then there exist three simple loops contained in the closure of $D$
and dual to the three edges incident to $v$.
Up to a little isotopy, these loops can be seen as loops in $\pi_1(\Sigma_0,v)$.
Up to orientation, each of them is the composition of the other two,
so at least one of them is orientation preserving in $M$.
This loop is orientation-preserving in $\Sigma$ and in $M$, 
and intersects $S(P)$ once: 
it easily follows that it bounds a disc in the ball $M\setminus
P$, which is absurd (since $P$ is minimal).
\end{proof}

\begin{lemma}\label{no_loops_in_p(G):lem}
Each edge in $p(G)$ has different endpoints.
\end{lemma}
\begin{proof}
Suppose that there exists an edge $e$ of $p(G)$ which joins a vertex $v$ to itself.
We have two cases depending on whether $v$ is 3-valent or 4-valent.
Suppose first that $v$ is 3-valent.
Since $\Sigma$ is orientable, the regular neighborhood ${\cal
R}_{\Sigma}(e)$ of $e$ in $\Sigma$ is an annulus.
Now there are two boundary components of ${\cal R}_{\Sigma}(e)$ in
$\Sigma$; one of these two loops does not intersect $S(P)$, so it is contained in a face of $P$.
Then there exists a face of $P$ incident to one vertex only: 
this contradicts the minimality of $P$.

We are left to deal with the case where $v$ is 4-valent.
We have two cases depending on whether the two germs of $e$ near $v$
lay on opposite sides with respect to $v$ or not.
If they do, the edge $e$ is the boundary of a face 
(not contained in $\Sigma$) incident to one vertex only: this contradicts the 
minimality of 
$P$.
In the second case, we
cannot choose a transverse orientation for ${\cal R}_{\Sigma}(e)$, 
because $P$ near $v$ lays (locally) on both sides of ${\cal R}_{\Sigma}(e)$.
Now, since $\Sigma$ is orientable, the regular neighborhood 
${\calR}_{\Sigma}(e)$ of $e$ in $\Sigma$ is an annulus.
Hence there are two boundary components of ${\cal R}_{\Sigma}(e)$ in
$\Sigma$; one of these two loops does not intersect $S(P)$, so it is contained in a face of $P$.
This loop is orientation reversing and bounds a disc: a contradiction.
\end{proof}

\begin{lemma}\label{external_length_4:lem}
If a connected component of $P\setminus\Sigma$ is a disc (so it is a
face of $P$ incident to $4$-valent vertices of $p(G)$ only), then it is incident
to at least $4$ vertices of $p(G)$ (with multiplicity).
\end{lemma}
\begin{proof}
If the disc is incident to 3 vertices at most (with multiplicity),
it is embedded by Lemma~\ref{no_loops_in_p(G):lem}, contradicting the
minimality of $P$.
\end{proof}

Now we are able to prove the four lemmas of Section~\ref{complexity:section}.
From now on, we suppose that $P$ has at most $6$ vertices.

\subsection{Proof of Lemma~\ref{general_sigma:lem}}

Recall that we want to prove that the Stiefel-Whitney surface $\Sigma$ contains at most 2 connected components and that $M$ has a minimal
standard spine with a Stiefel-Whitney surface having Euler characteristic equal to zero.
We will first suppose that $\Sigma$ is not connected, proving that
there are at most 2 components, and then we will prove that, up to
changing $P$, the Euler characteristic of $\Sigma$ is zero.

So let us suppose that $\Sigma$ is not connected.
Note that each component of $p(G)$ contains an even number of 3-valent vertices.
Hence, by Lemmas~\ref{one:vertex:lem}
and~\ref{existence_3_valent_vertices:lem}, 
each component of $p(G)$ contained in a component $\Sigma_0$ of 
$\Sigma$ contains at least one
4-valent vertex and a pair of 3-valent vertices; so $\Sigma_0$
contains at least 3 vertices of $P$.
Since $P$ has $6$ vertices at most, $\Sigma$ has two connected
components, each containing exactly 3 vertices of $P$.

Now, let us consider the Euler characteristic of $\Sigma$.

\paragraph{If $\Sigma$ has two components}
Let us concentrate on a connected component $\Sigma_0$ of $\Sigma$.
The Euler characteristic $\chi(\Sigma_0)$ can be computed using the
cellularization induced on $\Sigma_0$ by $P$.
The number of vertices is 3; so, since there are one 4-valent and two 3-valent vertices, then the number of
edges of $S(P)\cap \Sigma_0$ is equal to 5.
Now, $3-5=-2$ and there is at least one disc, so $\chi(\Sigma_0)\geqslant -1$.
We have already noted that each component of $\Sigma$ is different
from the sphere; 
so, since $\Sigma_0$ is orientable, then $\Sigma_0$ is a torus.

\paragraph{If $\Sigma$ is connected}
Let $g$ be the genus of the connected surface $\Sigma$.
We have already noted that $\Sigma$ is different from the sphere.
Let us suppose that $\Sigma$ is not a torus (\emph{i.e.}~$g\geqslant 2$).
We will first prove that $g < 4$, and then we will prove that 
the two remaining cases ($g=2,3$) are forbidden.
Let $v_3$ be the number of pairs of 3-valent vertices and $v_4$
the number of 4-valent vertices of $p(G)$.
As above, $\chi(\Sigma)$ can be computed using the cellularization
induced on $\Sigma$ by $P$.
The number of vertices is $2v_3+v_4$, thus we have
\begin{equation} \label{prima:eqn}
2v_3+v_4 \leqslant 6,
\end{equation}
where equality holds when all vertices of $P$ lie in $\Sigma$.
Since there are $v_4$ four-valent and $2v_3$ tri-valent vertices, the
number of edges of $S(P)\cap \Sigma$ is equal to 
$\frac{3(2v_3)+4v_4}{2}=3v_3+2v_4$.
Thus we have $\chi(\Sigma) = (2v_3+v_4)-(3v_3+2v_4)+f = f-v_3-v_4$,
where $f$ is the number of discs in $\Sigma\setminus S(P)$, so
\begin{equation} \label{seconda:eqn}
v_3+v_4 = 2(g-1)+f.
\end{equation}

The number of vertices of $P$ is greater than or equal to $v_3+v_4$, and
if $g\geqslant 4$, then $v_3+v_4\geqslant 6+f>6$, a contradiction.
So we are left to deal with a surface $\Sigma$ of genus 2 or 3.

\paragraph{If $\Sigma$ has genus 3}
If there is at least a 3-valent vertex ($v_3>0$), then $f\geqslant 2$ by
Lemma~\ref{3_valent_two_discs:lem}, so $v_3+v_4 \geqslant 6$ by~(\ref{seconda:eqn}).
Hence $2v_3+v_4>v_3+v_4 \geqslant 6$, contradicting~(\ref{prima:eqn}).
Therefore there are only 4-valent vertices ($v_3=0$), which implies
that $S(P)=p(G)$ and $P\setminus\Sigma$ consists of faces.
Since $\chi(P)=1$ and $\chi(\Sigma)=-4$, there are $5$ faces in $P\setminus\Sigma$.
Each (4-valent) vertex (of $p(G)$) is adjacent to exactly 2 germs of faces
of $P\setminus\Sigma$.
By Lemma~\ref{external_length_4:lem}, there should be at least
$5\cdot 4=20$ germs of such faces; but there are at most $6$ vertices in
$p(G)$, so there are at most $12$ germs of faces of $P\setminus\Sigma$.
A contradiction.

\paragraph{If $\Sigma$ has genus 2}
By Lemma~\ref{one:vertex:lem} we have $v_4>0$, so $v_3$ may be equal
to 0, 1, or 2 by~($\ref{prima:eqn}$).

\subparagraph{Case $v_3=2$}
We have $f\geqslant 2$ by Lemma~\ref{3_valent_two_discs:lem} and
$v_4=f$ by~(\ref{seconda:eqn}). Then~(\ref{prima:eqn})
implies that $v_4=2$. Thus $P$ has $2+4=6$ vertices and $7$ faces
(since $\chi(P)=1$), two of them in $\Sigma$ and 5 in $P\setminus\Sigma$.
These 5 faces of $P\setminus\Sigma$ may be incident three times to
each 3-valent vertex of $p(G)$ and twice to each 4-valent vertex
of $p(G)$.
Summing up, we obtain $16$ vertices (with multiplicity) 
to which the 5 faces are incident; so, among them, there exists a face 
incident to at most 3 vertices.
Such a face is embedded by Lemma~\ref{no_loops_in_p(G):lem}, in contrast to
the minimality of $P$.

\subparagraph{Case $v_3=1$}
We have $f\geqslant 2$ by Lemma~\ref{3_valent_two_discs:lem},
so $v_3+v_4\geqslant 4$ by~(\ref{seconda:eqn}).
Now there are two cases, depending on $v_4$.

If $v_4=4$, then $f=3$ and all $6$ vertices of $P$ belongs to $\Sigma$.
Since $\chi(P)=1$, there are $7$ faces in $P$, four of them in $P\setminus\Sigma$.
These $4$ faces are incident to $14$ vertices of $p(G)$ (with
multiplicity), so there exists a face incident to at most 3 
(which is embedded by Lemma~\ref{no_loops_in_p(G):lem}), a contradiction.

If $v_4=3$, then $f=2$ and $P$ has $5$ or $6$ vertices. If it has $5$
vertices (all contained in $\Sigma$), it has $6$ faces, $4$ of
them lying outside $\Sigma$.
These 4 faces are incident to $2\cdot v_4+3\cdot 2v_3 = 12$ vertices (with multiplicity), thus
there exists a face incident to at most 3 vertices, a contradiction.
If $P$ has $6$ vertices (one of them outside $\Sigma$), it has $7$
faces, $5$ of them lying outside $\Sigma$. These $5$ faces are
incident to $2\cdot v_4+3\cdot 2v_3+6 = 18$ vertices of $P$ (with multiplicity, the
vertex outside $\Sigma$ being counted 6 times), so there is a face
incident to at most 3, a contradiction.

\subparagraph{Case $v_3=0$}
We have $f=v_4-2$. There are only 4-valent vertices, then $S(P)=p(G)$ and
$P\setminus\Sigma$ consists of 3 disjoint discs (since $\chi(P)=1$ and $\chi(\Sigma)=-2$).
By Lemma~\ref{external_length_4:lem}, these discs are incident to at least $3\cdot
4=12$ vertices (with multiplicity), so $v_4=6$.
Now, let us consider the surface $\partial\calR(\Sigma)$, which is a double covering of $\Sigma$.
There are two cases depending on whether $\partial\calR(\Sigma)$ is connected or not.

Suppose first that $\partial\calR(\Sigma)$ is not connected, so it has two components which have genus 2.
Note now that $G$ is the disjoint union of three circles.
This fact contradicts Lemma~\ref{sphere_minus_discs:lem}, because two
surfaces of genus 2 minus three circles cannot be embedded in a
sphere.

\begin{figure}
\begin{center}
\mettifig{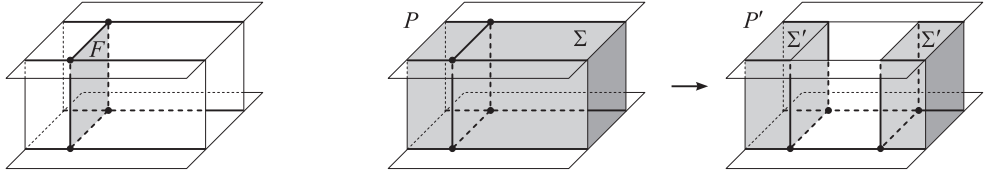}
\caption{If a face $F$ not contained in $\Sigma$ is embedded, we can
modify $P$ to a $P'$ with $\chi(\Sigma') = \chi(\Sigma)+2$.}
\label{embedded_disc:fig}
\end{center}
\end{figure}

\begin{figure}
\begin{center}
\mettifig{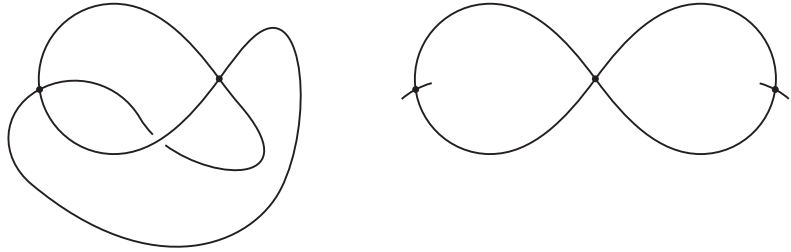}
\caption{The two cases for the boundary of each face not contained in $\Sigma$ (the dots are vertices of $P$).}
\label{two_cases_boundary_region:fig}
\end{center}
\end{figure}

\begin{figure}
\begin{center}
\mettifig{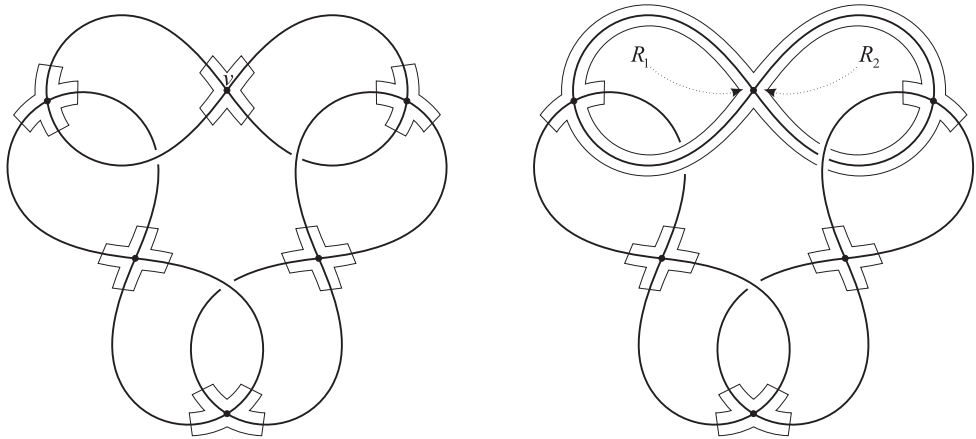}
\caption{The singular set of $P$ if $\Sigma$ contains only 4-valent
vertices (left).
Since $\Sigma$ is orientable, we can fix some matchings along the
edges (right).} 
\label{sing_set_4_val:fig}
\end{center}
\end{figure}

Suppose now that $\partial\calR(\Sigma)$ is connected, so it has genus 3.
By Lemma~\ref{external_length_4:lem}, each face which is not contained
in $\Sigma$ must be incident to at least four vertices.
Since each of the six vertices is adjacent to two faces (possibly the same)
not contained in $\Sigma$, 
then each of the three faces not contained in
$\Sigma$ is incident to four vertices.
Let us suppose first that there exists a face not contained in
$\Sigma$ 
which is embedded.
In such a case, by applying the move shown in
Fig.~\ref{embedded_disc:fig}, we obtain a 
spine $P'$ of $M$, with the same number of vertices
of $P$, with a new surface $\Sigma'$ which is a torus, and we are done.
So we are left to deal with the last case: namely, we suppose that all
faces not contained in $\Sigma$ are not embedded.
Lemma~\ref{no_loops_in_p(G):lem} easily implies that, for the boundary of
each face 
not contained in $\Sigma$, we have one the two cases 
shown in Fig.~\ref{two_cases_boundary_region:fig}.
Since $p(G)=S(P)$ is connected, the first case is not possible, so we are left to deal with the second one and $p(G)=S(P)$ appears as in 
Fig.~\ref{sing_set_4_val:fig}-left (we have shown also the neighborhood of the vertices in $\Sigma$).

Since $P$ is standard, then, to define $P$ uniquely, it is enough to say how 
the neighborhoods of the vertices match to each other along the
edges.
Since $\Sigma$ is orientable, we can suppose (up to symmetry) that, 
along the edges incident to the vertex $v$, the matchings are those shown
in Fig.~\ref{sing_set_4_val:fig}-right.
Now, note that all four faces contained in $\Sigma$ are incident to 24 vertices 
(with multiplicity).
For the two faces $R_1$ and $R_2$, indicated in Fig.~\ref{sing_set_4_val:fig}-right, 
we have two cases.
\begin{description}

\item{If $R_1 = R_2$.} The face $R_1=R_2$ is incident to at least 14 vertices, 
then the other 3 faces contained in $\Sigma$ are incident to
at most $24-14=10$ vertices: a contradiction.

\item{If $R_1 \neq R_2$.} Each of the faces $R_1$ and $R_2$ are incident to at least 10 vertices, so the other 2 faces contained in
$\Sigma$ are incident to at most $24-20=4$ vertices: a contradiction.

\end{description}

\subsection{Proof of Lemma~\ref{2_tori_sigma:lem}}

Recall that we want to prove that, if $\Sigma$ consists of two tori, $\calR(\Sigma)$ 
consists of two copies of $T\timtil I$ and each of the components of $\calR_P(\Sigma)$
appears as shown in Fig.~\ref{nhbd_T_non_orient:fig}.
It has been shown at the beginning of the proof of
Lemma~\ref{general_sigma:lem} that each component of $\Sigma$ contains
3 vertices of $P$.
As said above, there are two interval bundles on the torus up to
homeomorphism, namely the orientable $T\times I$ 
and the non-orientable $T\timtil I$.
If a component $\Sigma_0$ of $\Sigma$ has an orientable
neighborhood, $\partial\calR(\Sigma_0)$ consists of two tori, each
containing a component of $G\cap \partial\calR(\Sigma_0)$ with at least two vertices
by Lemma~\ref{discs:lem}: thus there are at least 4 vertices in
$\Sigma_0$, a contradiction.
Therefore each component of $\Sigma$ has a non-orientable neighborhood.

Let $\Sigma_0$ be a component of $\Sigma$. Since $\Sigma_0$
 contains 3
vertices, Lemma~\ref{discs:lem} implies that $G$ is a
$\theta$-graph in the torus $\partial\calR(\Sigma_0)$, and that $p$
maps the $\theta$-graph in $\Sigma_0$ producing one 4-valent vertex.
It is now easy to show that $\calR_P(\Sigma_0)$
 appears as shown in 
Fig.~\ref{nhbd_T_non_orient:fig}, so we leave it to the reader.

\subsection{Proof of Lemma~\ref{product_sigma:lem}}

Recall that we are analyzing the case when $\Sigma$ is one torus and 
$\calR(\Sigma)$ is the orientable $T\times I$, and we want to prove that
$\calR_P(\Sigma)$ appears as in Fig~\ref{nhbd_T_orient:fig}.
Lemma~\ref{one:vertex:lem} implies that $G$ contains at most
$6-1$ vertices, hence it contains 0, 2 or 4 vertices (being 3-valent).
Now, since $\calR(\Sigma)$ is orientable $G$ has two components.
It follows from Lemma~\ref{discs:lem} that $G$ consists of two
$\theta$-graphs, each mapped injectively into a $\theta$-graph in
$\Sigma$. Two $\theta$-graphs in a torus intersect transversely in at
least two points, and they intersect in exactly two only if they share
two slopes, \emph{i.e.}~if they are either isotopic or related by a flip.
Therefore $\calR_P(\Sigma)$ is one of the polyhedra shown in
Fig.~\ref{nhbd_T_orient:fig}.

\subsection{Proof of Lemma~\ref{twisted_sigma:lem}}

Recall that we are analyzing the case when $\Sigma$ is one torus and
$\calR(\Sigma)$ is the non-orientable $T\timtil I$, and we want
to prove that $M$ has a minimal standard spine with a Stiefel-Whitney 
surface $\Sigma$ such that $\calR(\Sigma) = T\timtil I$ and
$\calR_P(\Sigma)$ is as shown in Fig.~\ref{nhbd_T_non_orient:fig}.
As above, Lemma~\ref{one:vertex:lem} implies that $G$ contains at most
$6-1$ vertices, hence it contains 0, 2 or 4 vertices (being 3-valent).
It follows from Lemma~\ref{discs:lem} that $G$ contains
2 or 4 vertices.
If $G$ contains 2 vertices, then it is a $\theta$-graph in the
torus $\partial\calR(\Sigma)$. Therefore $M$ is obtained assembling
$\calR(\Sigma)$ and $N=M\setminus \calR(\Sigma)$, each manifold having
one torus boundary component marked with $G$.
Moreover, $\calR_P(\Sigma)$ and 
$Q=\closure{P\setminus\calR_P(\Sigma)}$ are skeleta for
$\calR(\Sigma)$ and $N$. Since $N$ is not a solid torus, Proposition~\ref{only:3:prop}
shows that $c(N)\geqslant 3$, thus $Q$ contains at least 3 vertices,
hence $\Sigma$ contains at most 3 vertices. 
As in the proof of Lemma~\ref{2_tori_sigma:lem}, we deduce that
$\calR(\Sigma)$ 
appears as shown in Fig.~\ref{nhbd_T_non_orient:fig}.

\paragraph{If $G$ contains 4 vertices}
To conclude the proof, we show that if $G$ contains 4
vertices, 
then we can modify $P$ to another spine $P'$ of $M$ with the same number
of vertices of $P$, with $\Sigma'\subset P'$ being a torus again,
such that $\calR(\Sigma') = T\timtil I$ and $G'$
contains two vertices.
Then the conclusion follows from the discussion above.

By applying Lemma~\ref{discs:lem} we get that
$\partial\calR(\Sigma)\setminus G$ is made of two open discs (say $D$ and $D'$).
Let us denote by, respectively, $e(D)$ and $e(D')$ the number of their edges (with
multiplicity). We have $e(D)+e(D')=12$, and we suppose $e(D)\leqslant
e(D')$, so $e(D)\leqslant 6$. 
Consider the polyhedron $P\cup D$. Then
$M\setminus (P\cup D)$
consists of two balls, one of them lying inside $\calR(\Sigma)$.
For each edge $s$ of $\partial D$, define $f_s$ to be the face
of $P\cup D$ incident to $s$ and contained in
$\calR(\Sigma)$. 
If $s$ has distinct endpoints $q_0, q_1$, then $f_s$
is incident to 4 distinct vertices $q_0,q_1,p(q_0),p(q_1)$ of $P\cup
D$. Now $(P\cup D)\setminus f_s$ is another spine of $M$ with $4-e(D)$ vertices less
than $P$, hence $e(D)\geqslant 4$ (using that $D$ is
embedded if $e(D)\leqslant 4$).
If $e(D)=4$, the spine
$(P\cup D)\setminus f_s$ is standard and minimal, and the new $G'$ has two vertices only.

\begin{figure}
\begin{center}
\mettifig{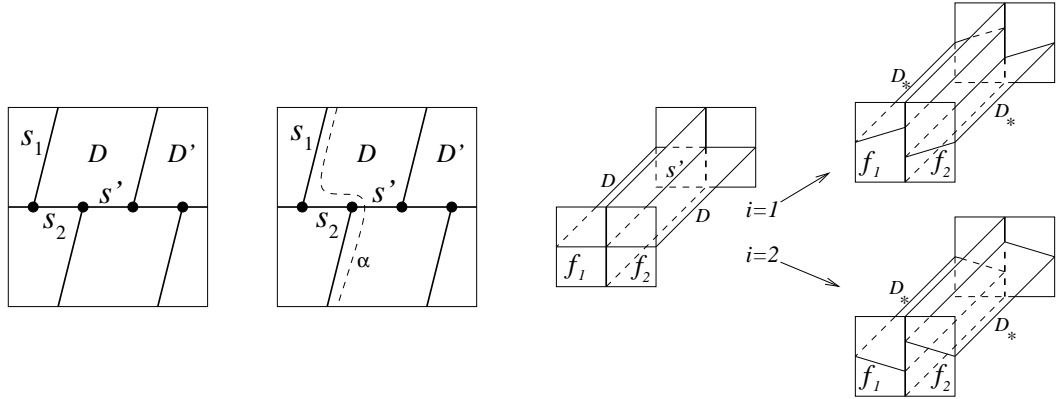, width = 14 cm}
\caption{The only possibility for $G$ with $e(D)>4$ (left), 
the loop $\alpha$ bounding a disc in $M\setminus P$ (centre),
and a perturbation of $D$ into $D_*$ so that 3 vertices of $\partial D_*$
are adjacent to $f_{e_i}$ (right).} \label{4vert:fig}
\end{center}
\end{figure}

There is only one case with $e(D)>4$, shown in Fig.~\ref{4vert:fig}-left (we have
$e(D)=e(D')=6$). Set $f_1 = f_{s_1}, f_2 = f_{s_2}$. Each
$f_i$ separates the two balls given by $M\setminus(P\cup D)$. 
and is incident to at least two vertices of
$\Sigma$. If each $f_i$ is incident only twice, then $p(s_1\cup
s_2)$ does not contain any 4-valent vertex.
The loop $\alpha\subset \partial\calR(\Sigma)$ shown in Fig.~\ref{4vert:fig}-centre then
projects to a simple loop $p(\alpha)$ in $\Sigma\setminus p(G)$
which bounds a disc in the ball $M\setminus P$ and
meets $p(G)$ in one point, which is absurd (since $P$ is minimal).
Therefore some $f_i$ is incident to at least 3
vertices of $p(G)$, for $i=1$ or~$2$.
The disc $D$ is not embedded, so we perturb it 
into an embedded $D_*$.
We can do it so that 3 vertices of $\partial D_*$ are adjacent to
$f_i$ (with a perturbation depending on $i$, see Fig.~\ref{4vert:fig}-right), 
thus $f_i$ is incident to at least 6 distinct vertices of
$P\cup D_*$.

If $f_i$ is incident to more than 6 vertices, then $(P\cup D_*)\setminus f_i$ is a
standard spine of $M$ with less vertices than $P$.
Therefore $f_i$ is incident to exactly 6 vertices and $(P\cup D_*)\setminus f_i$ 
is the required minimal standard spine of $M$ 
(with the same number of vertices of $P$), with a new $G'$ having 2 vertices.

{\small

\vspace{1.5 cm}

\noindent
\hspace*{6cm}Dipartimento di Matematica\\ 
\hspace*{6cm}Universit\`a di Pisa\\ 
\hspace*{6cm}Via F. Buonarroti 2\\ 
\hspace*{6cm}56127 Pisa, Italy\\ 
\hspace*{6cm}amendola@mail.dm.unipi.it\\
\hspace*{6cm}martelli@mail.dm.unipi.it
}

\end{document}